\documentclass[preprint,12pt,3p]{elsarticle}

\usepackage{amssymb}
\usepackage{amsmath,amsthm}
\usepackage{mathrsfs}
\usepackage{hyperref}
\usepackage{cleveref}
\usepackage[all,cmtip]{xy}
\usepackage{epsf, graphicx}
\usepackage{latexsym,amsfonts,amsbsy,amssymb}
\usepackage{geometry}
\usepackage{titletoc}
\usepackage{rotating}
\usepackage{algorithm}
\usepackage{algorithmic}
\usepackage{amscd}

\makeatletter

\@addtoreset{equation}{section} \makeatother
\newtheorem{theorem}{Theorem}

\newtheorem{conjecture}[theorem]{Conjecture}
\newtheorem{corollary}[theorem]{Corollary}
\newtheorem{remark}[theorem]{Remark}

\def\IBr{{\rm IBr}}

\def\O{\mathcal{O}}

\def\Z{\mathbb{Z}}
\def\Q{\mathbb{Q}}

\makeatletter
\def\ps@pprintTitle{%
\let\@oddhead\@empty
\let\@evenhead\@empty
\def\@oddfoot{\reset@font\hfil\thepage\hfil}
\let\@evenfoot\@oddfoot
}
\makeatother

\begin{document}

\begin{frontmatter}
	
	\title{Virtual Morita equivalences and Brauer character bijections}

	\author{Xin Huang}
	
	\begin{abstract}
	We extend a theorem of Kessar and Linckelmann concerning Morita equivalences and Brauer character bijections between blocks to virtual Morita equivalences. As a corollary, we obtain that Navarro's refinement of Alperin's weight conjecture holds for blocks with cyclic and Klein four defect groups, blocks of symmetric and alternating groups with abelian defect groups, and $p$-blocks of ${\rm SL}_2(q)$ and ${\rm GL}_2(q)$, where $p|q$.
	\end{abstract}
	

\end{frontmatter}

We first review some notation in \cite{Kessar_Linckelmann} which is necessary for stating the main result. Let $p$ be a prime number and $(K, \O, k)$ a
$p$-modular system; that is, $\O$ is a complete discrete valuation ring 
with residue field $k=\O/J(\O)$ of characteristic $p$ and field of 
fractions $K$ of characteristic zero. Assume that $k$ is a perfect field. We say that a positive integer $n$ is {\it large enough} for a finite group $G$ if $n$ is a multiple of the order of $G$. Let $\bar{K}$ be an algebraic closure of $K$. Let $n$ be a positive 
integer, denote by $\Q_n$ the  $n$-th cyclotomic subfield of $\bar{K}$, 
and let $k'$ be a splitting field of the polynomial $x^{n_{p'}}-1$ over 
$k$, where $n_{p'}$ denotes the $p'$-part of $n$. Denote by $\mathcal{H}_n$ the subgroup of ${\rm Gal}(\Q_n/\Q)$ consisting of 
those automorphisms $\alpha$ for which there exists a non-negative 
integer $u$ such that $\alpha(\delta) =\delta^{p^u}$ for all 
$n_{p'}$-roots of unity $\delta$ in $\Q_n$. Denote by $\mathcal{H}_{n,k}$ the 
subgroup of $\mathcal{H}_n$ consisting of those automorphisms $\alpha$ for 
which there exists a non-negative integer $u$ and an element  $\tau \in 
{\rm Gal}(k'/k)  $ such that $\alpha (\delta) =\delta^{p^u}$  for all 
$n_{p'}$-roots of unity $\delta$ in $\Q_n$  and $\tau( \eta) =\eta^{p^u}$  for all $n_{p'} $-roots of  unity  $\eta$ in $k'$.  Note that   $\mathcal{H}_{n,k} $ is independent of the choice of a splitting field  
$k'$ of $x^{n_{p'}}-1$ over $k$.

For a finite group $G$, denote by $\IBr(G)$ the set 
of irreducible Brauer characters of $G$ interpreted as functions from 
the set of $p$-regular elements of $G$ to $\bar K$.  If $b$ is a central 
idempotent of $\O G$, then we denote by $\IBr(G,b) $ the subset of 
$\IBr(G)$ consisting of the Brauer characters of simple $k'Gb $-modules
for  any sufficiently large field  $k'$ containing $k$.
The group $\mathcal{H}_{n,k}$ acts on $\IBr(G,b)$ via  
${}^\sigma\varphi(g):= \sigma(\varphi(g))$,  for any $\varphi \in \IBr(G,b) $, $\sigma\in\mathcal{H}_{n,k}$ and any $p$-regular element $g\in G$. A {\it $p$-weight} of $G$ is a pair $(Q,\phi)$, where $Q$ is a $p$-subgroup of $G$ and $\phi\in \IBr(N_G(Q))$ is such that the kernel of $\phi$ contains $Q$ and such that the $p$-parts of $\phi(1)$ and the order of $N_G(Q)/Q$ are the same. 
Denote by $\bar{b}$ the image of $b$ in $kG$.

We generalise the second statement in \cite[Theorem 3.4]{Kessar_Linckelmann} by considering virtual Morita equivalences instead of Morita equivalences. For the definition of a virtual Morita equivalence we refer to \cite[Definition 1.1]{Kessar_Linckelmann} or \cite[Definition 9.1.4]{Linckelmann}. Morita equivalences, Rickard equivalences and $p$-permutation equivalences all induce virtual Morita equivalences (see \cite[Remark 1.2]{Kessar_Linckelmann}).

\begin{theorem}\label{main}
Let $G$ and $H$ be finite groups and let $n$ be a positive integer large enough for $G$ and $H$. Let $b$ and $c$ be block idempotents of $\O G$ and $\O H$ respectively. Suppose that $\O$ is absolutely unramified (i.e. $J(\O)=p\O$). If there is a virtual Morita equivalence between $\O Gb$ and $\O Hc$ given by a virtual bimodule and its dual, then there is a bijection $I:\IBr(G,b)\to \IBr(H,c)$ satisfying ${}^\sigma I(\varphi)=I({}^\sigma\varphi)$ for all $\varphi\in \IBr(G,b)$ and $\sigma\in \mathcal{H}_{n,k}$. 
\end{theorem}

\noindent{\it Proof.} Without loss of generality we may assume that $k$ is finite (because every finite group has a finite splitting field). Let $(K',\O', k')$ be  an extension of  the $p$-modular system
$(K,\O,k)$ (see \cite[Definition 3.1]{Kessar_Linckelmann}) such that $K'\subseteq$ $\bar K$ and such that the 
extension $K'/K$ is normal. Suppose that  $k'$ is finite, and that 
$K'$  contains primitive  $n$-th roots of 
unity. Since $k'$ contains  enough roots of unity, we may identify $\IBr(G,b) $  
(resp. $\IBr(H,c)$) with the Brauer characters of simple 
$k'G\bar{b} $-modules (resp. $k'H\bar{c}$-modules). Let $\Gamma:={\rm Gal}(k'/k)$. Denote by $\mathcal{S}(k'G\bar{b})$ the set of isomorphism classes of simple $k'G\bar{b}$-modules. The group $\Gamma$ acts on the set $\mathcal{S}(k'G\bar{b})$ (see \cite[Notation 1.15]{Kessar_Linckelmann} for the action). By \cite[Lemma 3.3 (b) and Lemma 3.2]{Kessar_Linckelmann}, to prove the theorem, it suffices to show that $\mathcal{S}(k'G\bar{b})$ and $\mathcal{S}(k'H\bar{c})$ are isomorphic as $\Gamma$-sets.

Now suppose that $X\in\mathcal{P}(\O Gb,\O Hc)$ (see \cite[Definition 1.1]{Kessar_Linckelmann}) and its dual induce a virtual Morita  
equivalence  between $\O Gb $ and $\O Hc $. Then $X':=\O' \otimes_\O X$ 
and its dual induce a virtual Morita equivalence between $\O'Gb $ and $\O' Hc$; $\bar{X}:= k\otimes_{\O} X $ and $\bar{X}':=  k'\otimes_{k} \bar X$ (together with their duals)
induce virtual Morita equivalences between $kG\bar{b}$ and $kH\bar{c}$ and between $k'G\bar{b}$ 
and $k'H\bar{c}$ respectively.

By \cite[Corollary 9.3.4]{Linckelmann}, $X'$ and its dual induce a perfect isometry $\Phi:\mathcal{R}(K'Hc)\to \mathcal{R}(K'Gb)$, where the notation $\mathcal{R}(-)$ means the Grothendieck group (see \cite[Definition 1.8.6]{Linckelmann}). By \cite[Corollary 9.2.7]{Linckelmann}, $\Phi$ induces a unique $\Z$-module isomorphism $\bar{\Phi}:\mathcal{R}(k'H\bar{c})\cong \mathcal{R}(k'G\bar{b})$ such that $d_G\circ\Phi=\bar{\Phi}\circ d_H$, where $d_G$ and $d_H$ are decomposition maps (see \cite[Definition 4.17.2]{Linckelmann}). 
Actually it easy to see that $\bar{\Phi}$ is induced by tensoring with $\bar{X}'$. The action of $\Gamma$ on $\mathcal{S}(k'G\bar{b})$ (resp. $\mathcal{S}(k'H\bar{c})$) extends linearly to an action on $\mathcal{R}(k'G\bar{b})$ (resp. $\mathcal{R}(k'H\bar{c})$), so we can regard $\mathcal{R}(k'G\bar{b})$ and $\mathcal{R}(k'H\bar{c})$ as $\Z\Gamma$-modules. It is easy to check that $\bar{\Phi}$ is a $\Z\Gamma$-module isomorphism. 
Indeed, since $\bar{X}'\cong k'\otimes_k \bar{X}$ as $(k'G\bar{b},k'H\bar{c})$-bimodules, we have ${}^\gamma \bar{X}'\cong \bar{X}'$ for any $\gamma\in \Gamma$ (see the paragraph after \cite[Lemma 6.1]{Kessar_Linckelmann}). It follows that for any $\gamma\in\Gamma$ and any simple $k'H\bar{c}$-module $S$, we have
$${}^\gamma(\bar{\Phi}([S]))={}^\gamma[\bar{X}'\otimes_{k'H\bar{c}} S]=[{}^\gamma\bar{X}'\otimes_{k'H\bar{c}} {}^\gamma S]= [\bar{X}'\otimes_{k'H\bar{c}} {}^\gamma S]=\bar{\Phi}({}^\gamma[S]).$$
Since $\mathcal{S}(k'H\bar{c})$ is a $\Gamma$-stable $\Z$-basis of $\mathcal{R}(k'H\bar{c})$, $\bar{\Phi}(\mathcal{S}(k'H\bar{c}))$ is a $\Gamma$-stable $\Z$-basis of $\mathcal{R}(k'G\bar{b})$ as well. On the other hand, $\mathcal{S}(k'G\bar{b})$ is another $\Gamma$-stable $\Z$-basis of $\mathcal{R}(k'G\bar{b})$. 

Since the character of the $\Z \Gamma$-module $\mathcal{R}(k'G\bar{b})$ is independent of the choice of a $\Z$-basis, we see that for any $\gamma\in \Gamma$, the number of fixed points of $\gamma$ on  $\bar{\Phi}(\mathcal{S}(k'H\bar{c}))$ equals the number on $\mathcal{S}(k'G\bar{b})$. Since $\Gamma$ is a cyclic group, by \cite[Lemma 13.23]{Isaacs}, $\bar{\Phi}(\mathcal{S}(k'H\bar{c}))$ and $\mathcal{S}(k'G\bar{b})$ are isomorphic as $\Gamma$-sets. Hence $\mathcal{S}(k'H\bar{c})$ and $\mathcal{S}(k'G\bar{b})$ are isomorphic as $\Gamma$-sets.   $\hfill\square$

\medskip In the end of \cite{Navarro}, Navarro proposed a refinement of Alperin's weight conjecture. In \cite[Conjecture]{Turull}, Turull stated the conjecture for a finite group. For a single block, Navarro's conjecture should be the following. 

\begin{conjecture}\label{conj:Navarro}
Let $\mathbf{Q}_p$ be the field of $p$-adic numbers, $G$ a finite group, $b$ a $p$-block idempotent of $G$ (over a suitable subring of $\bar{\mathbf{Q}}_p$ which is large enough for $G$), and ${\rm Gal}(\bar{\mathbf{Q}}_p/\mathbf{Q}_p)_b$ the stabliser of $b$ in ${\rm Gal}(\bar{\mathbf{Q}}_p/\mathbf{Q}_p)$. Let ${\rm W}(G,b)$ be the set of $G$-conjugacy classes of $p$-weights of $G$ in $b$. Then there is a bijection $\IBr(G,b)\to {\rm W}(G,b)$ commuting with the action of ${\rm Gal}(\bar{\mathbf{Q}}_p/\mathbf{Q}_p)_b$.
\end{conjecture}

If Conjecture \ref{conj:Navarro} holds for all blocks of a finite group $G$ then \cite[Conjecture]{Turull} holds for $G$, and vice versa. 
Combined with some results on Brou\'{e}'s abelian defect group conjecture over arbitrary complete discrete valuation rings (\cite[Theorem 1.12]{Kessar_Linckelmann}, \cite[Theorem 3]{Huang}, \cite[Theorem 7.6]{CR}, \cite[Theorem 1.1]{H21} and \cite[Theorem 1.1]{HLZ}), Theorem \ref{main} implies the following new result.
\begin{corollary}\label{corollary}
Conjecture \ref{conj:Navarro} holds for 
\begin{enumerate}[\rm(i)]
	\item blocks with cyclic and Klein four defect groups;
	\item blocks of symmetric and alternating groups with abelian defect groups;
		\item $p$-blocks of ${\rm SL}_2(q)$ and ${\rm GL}_2(q)$, where $p|q$.
	\end{enumerate}
\end{corollary}

\begin{remark}
{\rm	When using Theorem \ref{main} to prove Conjecture \ref{conj:Navarro} for a block $b$ of a finite group $G$ over a large enough complete discrete valuation ring $\O'$, we choose the ring $\O$ in Theorem \ref{main} to be the minimal complete discrete valuation ring of $b$ in $\O'$ (see \cite[Definition 1.8]{Kessar_Linckelmann}). The  action of $\mathcal{H}_n$ on $\IBr(G)$ induces an action of 
	$\mathcal{H}_n $ on the set of blocks of $\O'G$.   Since $\O$ is a minimal 
	complete discrete valuation ring for $b$, $k$  is  a finite field  
	and consequently a splitting field of $x^{n_{p'}}-1 $ over $k $ is also 
	finite.  Let $|k|=p^d$. Then  $\mathcal{H}_{n,k}$ consists of  precisely 
	those elements $\alpha$ of ${\rm Gal}(\Q_n/\Q)$  for which there exists a 
	non-negative integer $u$ such that $\sigma(\delta)=$ $\delta^{p^{ud}}$  
	for all $n_{p'} $-roots of unity  in $\Q_n$. It follows that 
	$\mathcal{H}_{n,k}$ is exactly the stabiliser of  $b$ in $\mathcal{H}_n$. For better understanding of the connection between Theorem \ref{main} and Conjecture \ref{conj:Navarro}, we recommend \cite[\S2]{Turull}, where we can learn more equivalent reformulations of Conjecture \ref{conj:Navarro}.}
\end{remark}

\medskip\noindent\textbf{Acknowledgements.} The author thanks Prof. Radha Kessar for encouraging comments and Prof. Markus Linckelmann for suggesting to state out definitely the block version of Navarro's refinement of Alperin's weight conjecture.


\begin{thebibliography}{99}
\addcontentsline{toc}{section}{\protect References}




\bibitem{CR}
J. Chuang, R. Rouquier, Derived equivalences for symmetric groups and $\mathfrak{sl}_2$-categorification, Ann. Math. {\bf 167} (2008) 245-298.

\bibitem{Huang}
X. Huang, Descent of equivalences for blocks with Klein four defect groups, J. Algebra {\bf 398} (2023) 898-905.

\bibitem{H21}
X. Huang, Descent of splendid Rickard equivalences in alternating groups, arXiv:2111.10922. 

\bibitem{HLZ}
X. Huang, P. Li, J. Zhang, The strengthened Brou\'{e} abelian defect group conjecture for ${\rm SL}(2,p^n)$ and ${\rm GL}(2,p^n)$, J. Algebra {\bf 633} (2023) 114-137.
 

\bibitem{Isaacs}
I. M. Isaacs, Character Theory of Finite Groups,  Academic Press, New York, 1976.

\bibitem{Kessar_Linckelmann}

R. Kessar, M. Linckelmann, Descent of equivalences and character bijections, Geometric and topological aspects of the representation theory of finite groups, Springer Proc. Math. Stat., {\bf 242} (2018) 181-212.

\bibitem{Linckelmann}
M. Linckelmann, The Block Theory of Finite Group Algebras I/II, London Math. Soc. Student Texts, vol. {\bf 91/92}, Cambridge University Press, 2018.

\bibitem{Navarro}
G. Navarro, The McKay conjecture and Galois automorphisms, Ann. Math. {\bf160} (2004) 1129-1140.

\bibitem{Turull}
A. Turull, The strengthened Alperin weight conjecture for $p$-solvable groups, J. Algebra {\bf 398} (2014) 469-480.  

\end{thebibliography}
\end{document}